\documentclass[11pt]{amsart}
\usepackage{mathtools, amssymb, amsfonts, amsthm, enumitem}
\usepackage{epsfig, psfrag, color, multicol, graphicx, tikz,pdfsync}
\usepackage{amsfonts}
\usepackage{epstopdf}
\usepackage{xspace}

\usepackage{hyperref}

\usepackage[all]{xy}
\usepackage[font=small,labelfont=bf]{caption}
\usepackage{hyperref}
\hypersetup{colorlinks,linkcolor={red},citecolor={olive},urlcolor={red}}

\usetikzlibrary{calc}
\mathtoolsset{showonlyrefs}

\pdfoutput=1

\def\.{\hskip.06cm}
\def\ts{\hskip.03cm}

\def\sq{\square}

\def\zz{\mathbb Z}
\def\nn{\mathbb N}

\def\rr{\mathbb R}

\def\wb{}

\def\la{\lambda}

\def\de{\delta}

\def\ve{\varepsilon}

\def\cS{\mathcal S}

\def\T{\mathbf{T}}

\def\<{\langle}
\def\>{\rangle}

\def\rM{ {\text {\rm {T}}}}
\def\rI{ {\text {\rm {G}}}}

\def\rT{{\text {\rm T} } }

\def\0{{\mathbf 0}}

\def\rTab{{\text{\rm T}}(\mathbf a, \mathbf b)}

\def\.{\hskip.06cm}
\def\ts{\hskip.03cm}

\def\bba{{\text{\bf a}}}
\def\bbb{{\text{\bf b}}}

\def\P{{\textup{\textsf{P}}}}

\def\nin{\noindent}

\def\limsup{\mathop{\rm lim\,sup}\limits}

\def\R{\mathbb{R}}

\def\P{\mathbb{P}}

\def\T{\mathcal{T}}

\def\a{\mathbf{a}}
\def\b{\mathbf{b}}

\def\r{\textbf{\textit{a}}}
\def\c{\textbf{\textit{b}}}

\def\M{\mathcal{T}}

\makeatletter
\newcommand*\rel@kern[1]{\kern#1\dimexpr\macc@kerna}
\newcommand*\widebar[1]{%
	\begingroup
	\def\mathaccent##1##2{%
		\rel@kern{0.8}%
		\overline{\rel@kern{-0.8}\macc@nucleus\rel@kern{0.2}}%
		\rel@kern{-0.2}%
	}%
	\macc@depth\@ne
	\let\math@bgroup\@empty \let\math@egroup\macc@set@skewchar
	\mathsurround\z@ \frozen@everymath{\mathgroup\macc@group\relax}%
	\macc@set@skewchar\relax
	\let\mathaccentV\macc@nested@a
	\macc@nested@a\relax111{#1}%
	\endgroup
}

\DeclareMathOperator*{\argmax}{arg\,max}

\usepackage{theoremref}

\newtheorem{theorem}{Theorem}
\numberwithin{theorem}{section}
\numberwithin{equation}{section}

\newtheorem{lemma}[theorem]{Lemma}
\newtheorem{prop}[theorem]{Proposition}
\newtheorem{corollary}[theorem]{Corollary}

\theoremstyle{definition}
\newtheorem{definition}[theorem]{Definition}

\allowdisplaybreaks

\begin{document}
	
\title[On the number of contingency tables]{On the number of contingency tables \\ and the independence heuristic}

\author[Hanbaek Lyu and Igor Pak]{Hanbaek Lyu$^\ast$ and Igor Pak$^\ast$}

	
\thanks{\thinspace \ \today}
\thanks{\thinspace ${\hspace{-.45ex}}^\star$Department of Mathematics,
UCLA, Los Angeles, CA, 90095;  \  Email:
\hskip.06cm  \texttt{\{hlyu,\ts pak\}@math.ucla.edu}\ts{}}


\begin{abstract}
We obtain sharp asymptotic estimates on the number of $n \times n$
contingency tables with two linear margins $Cn$ and~$BCn$.  The results
imply a second order phase transition on the number of such contingency
tables, with a critical value at \ts $B_{c}:=1 + \sqrt{1+1/C}$.
As a consequence, for \ts $B>B_{c}$, we prove that the classical
\emph{independence heuristic} leads to a large undercounting.
\end{abstract}
	
	${}$
	\vspace{-1.70cm}
	${}$
\maketitle

\section{Introduction}\label{Introduction}

Sometimes a conjecture is more than a straightforward claim to be proved
or disproved.  A conjecture can also represent an invitation to understand
a certain phenomenon, a challenge to be confirmed or refuted in every
particular instance.  Regardless of whether such a conjecture is true or false,
the advances toward resolution can often reveal the underlying nature
of the objects.

This paper concerns with the \emph{independence heuristic} for
approximating the number of contingency tables, introduced by
I.~J.~Good as far back as in~1950.  The independence heuristic
has been both proved and disproved in several extreme cases.
This paper investigates an intermediate case of the margins when
the asymptotics are very subtle.  Unreachable until now
with the existing techniques, the results are quite surprising,
providing a new piece of the puzzle.

\smallskip

Let \ts $\r=(a_1,\ldots,a_m) \in \nn^{m}$ \ts and \ts
$\c=(b_1,\ldots,b_n)\in \nn^{n}$, such that
\ts $a_1+\ldots + a_m = b_1+\ldots +b_n=N$.
A \emph{contingency table} with \emph{margins}~$(\r,\c)$ is a
$m\times n$ matrix \ts $X=\bigl(x_{ij}\bigr)$, s.t.\
$x_{ij}\in \nn$, and
\begin{align}\label{eq:def_contingency_set}
\sum_{j=1}^{n} \. x_{ij} \, = \, a_i\,, \quad \
\sum_{i=1}^{m} \. x_{ij} \, = \, b_j\, \quad \text{for all} \quad 1\le i\le n, \, \ 1\le j \le m\ts.
\end{align}
Denote by \ts $\M(\r,\c)$ \ts the set of such tables, and let \ts
$\rT(\r,\c) = \bigl|\M(\r,\c)\bigr|$.

Computing and approximating $\rT(\r,\c)$
is a fundamental in Statistics and a classical problem in Combinatorics,
with many connections and applications to other fields,
see e.g.~\cite{DG} (see also~\cite{BLP} for recent references).
While there are a number of algorithmic approaches and asymptotic
results for small margins, the lower and upper bounds for
large margins remain far apart, see~\cite{BLP}. In fact, there
is a dearth of asymptotic tools in the latter case, and very
little hope to get a tight asymptotic bound in full generality.

The \emph{independence heuristic} is a classical approximation formula:
\begin{equation}\label{eq:IH}
\rM(\r,\c) \. \approx \. \rI(\r,\c),
\end{equation}
where
\begin{equation}\label{eq:IH-def}
\rI(\r,\c)\, := \, \binom{N+mn-1}{mn-1}^{-1} \, \prod_{i=1}^{m}  \binom{a_{i}+n-1}{n-1} \,
\prod_{j=1}^{n} \binom{b_{j}+m-1}{m-1}\ts.
\end{equation}
The idea behind the independence heuristic is the asymptotic
independence of rows and columns of random
continency tables $X \in \M(\r,\c)$, see~$\S$\ref{ss:main-reason}.
We postpone the history of~\eqref{eq:IH} and numerical examples
until~$\S$\ref{ss:finrem-hist}.

For the uniform margins, the independence heuristic
was studied by Canfield and McKay \cite{CM}.  In particular,
for $m=n$, $a_i = b_i=Bn$,\footnote{To simplify the presentation,
throughout the introduction we drop the floor/ceiling notation,
and use $Bn$, $\ve\ts n$, $n^\de$, etc., to mean the
nearest integer to these values. } they prove that
\smallskip
\begin{equation}\label{eq:CM}
\rT(\r,\c) \. \sim \. \sqrt{e} \.\cdot\. \rI(\r,\c) \ \ \ \text{as} \ \ n \to \infty,
\end{equation}
where \ts 
$B>0$ is a fixed constant.  The same asymptotics~\eqref{eq:CM} was
proved by Greenhill and McKay \cite{GM} for small margins: \ts
$\max\{a_i\}\cdot \max\{b_j\} = o\bigl(N^{2/3}\bigr)$.

In the opposite direction, Barvinok proved that the independence
heuristic~\eqref{eq:IH} fails for nonuniform ``cloned margins''.
In a notable special case, for $m=n$,
$\r=\c=(Bn,\ldots,Bn,n,\ldots,n)$, with~$\ve \ts n$ and $(1-\ve)\ts n$
rows/columns of each sum, $B>1$, he proves:
$$
\lim_{n\to \infty} \. \frac{1}{n^2} \. \log \rT(\r,\c) \, > \, \lim_{n\to \infty} \. \frac{1}{n^2} \. \log \rI(\r,\c).
$$
In other words, the independence heuristic greatly undercounts
the number of contingency tables for constant fraction of each sum.

In this paper we consider an intermediate case $m=n$,
$\r=\c=(Bn,\ldots,Bn,n,\ldots,n)$, where~$n^\de$ rows/columns
have larger sums, and fixed $B>1$, $0\le \de<1$. In this case it
is known and easy to see that
$$
\lim_{n\to \infty} \. \frac{1}{n^2} \. \log \rT(\r,\c) \, =
\, \lim_{n\to \infty} \. \frac{1}{n^2} \. \log \rI(\r,\c)\, =
\, 2\. \log 2.
$$
We show that the independence heuristic works fairly well in this case
as the second terms of the asymptotics have the same order:
$$
\aligned
\log \rT(\r,\c) \, & = \, (2\.\log 2) \. n^2 \, + \, \Theta\bigl(n^{1+\de}\bigr), \\
\log \rI(\r,\c) \, & = \, (2\. \log 2) \. n^2 \, + \, \Theta\bigl(n^{1+\de}\bigr).
\endaligned
$$
But the similarities stop when we compute the exact constant implied by
the $\Theta(\cdot)$ notation.  We present the detailed results in the
next section, but here is the qualitative version of (the special case of)
the main theorem.

\smallskip

\begin{corollary}  \label{cor:main}
Fix $0<\de<1$, $B>1$, and denote \ts $B_c:=1+\sqrt{2}$.
Let $m=n$, $\r=\c=(Bn,\ldots,Bn,n,\ldots,n)$ with $n^\de$ sums $Bn$.  Then:
$$
\lim_{n\rightarrow\infty}\. \frac{1}{n^{1+\delta}} \, \log \frac{\rT(\r,\c)}{\rI(\r\,\c)}  \quad
\begin{cases}
\  =	\, 0 \ & \text{if \ \. $1<B \le B_{c}$} \\
\   > \, 0 \ & \text{if \ \. $B>B_{c}$}
\end{cases}
$$
\end{corollary}

\smallskip

This is quite surprising since the independence heuristic does not ``notice'' the phase
transition at~$B_c$ and changes smoothly with~$B$.  The corollary then implies the second
order phase transitions for the number of contingency tables, see the discussion below.

The significance of the critical value \ts $B_c=1+\sqrt{2}$ \ts for the distribution of
random contingency tables has already been predicted in~\cite{B3} and proved in~\cite{DLP19},
but until now they never appeared in the context of counting contingency tables.

To summarize the idea of the proof, we combined Barvinok's classical bounds and our
previous results on the distribution of entries in contingency tables.  We then use
\emph{self-reduction} to derive the asymptotics for the number of contingency tables.
Put succinctly, the difference in these distributions before and after the
phase transition then amplifies the undercounting by the independence heuristic.

\bigskip

\section{Main results}\label{sec:main}

\subsection{Barvinok margins}\label{ss:main-thm}
Fix parameters \ts $0\le \delta \le 1$, $B\ge 1$ and $C>0$. As in the introduction,
define \emph{Barvinok margins}
\begin{align}\label{eq:def_barvinok_margin}
\r\. = \. \c \. := \. \bigl(\lfloor BCn \rfloor ,\ldots,\lfloor BCn\rfloor ,\lfloor Cn \rfloor, \ldots, \lfloor Cn \rfloor \bigr)
\, \in \, \nn^{\lfloor n^{\delta} \rfloor \ts +\ts n}\ts.
\end{align}
To simplify the notation, for the Barvinok margins we write \ts $\M_{n,\delta}(B,C)$,
\ts $\rT_{n,\delta}(B,C)$ \ts and \ts $\rI_{n,\delta}(B,C)$.  Formally, \ts $\M_{n,\delta}(B,C)$ \ts
is the set of contingency tables whose first \ts $\lfloor n^{\delta} \rfloor$ \ts rows
and columns have sums \ts $\lfloor BCn \rfloor$, and the other $n$~rows and columns
have sums \ts $\lfloor Cn \rfloor$.  Similarly, \ts $\rT_{n,\delta}(B,C)= |\M_{n,\delta}(B,C)|$,
and \ts $\rI_{n,\delta}(B,C)$ \ts is the corresponding independence heuristic
approximation~\eqref{eq:IH-def}.
	
\smallskip

The main result of this paper is a sharp asymptotics for the number
\ts $\rT_{n,\delta}(B,C)$ \ts of contingency tables for Barvinok's margins.
The result establishes a phase transition at a critical value \ts $B_{c}=1+\sqrt{1+1/C}$,
where the second order term in \ts $\log \rT_{n,\delta}(B,C)$ \ts
grows in~$B$ for $B<B_{c}$, but remains constant for \ts $B>B_{c}$.

\smallskip	
	
\begin{theorem}[Main theorem]\label{thm:main_volume}
Fix $0< \delta<1$, $B,C>0$, and $n\ge 1$. Let \ts $B_{c}=1+\sqrt{1+1/C}$ \ts
and denote \ts $f(x) := (x+1)\log(x+1) - x\log x$. 

\smallskip
\nin
$(i)$ \. For \ts $B \le B_{c}$, we have:
\begin{align}
\log \rT_{n,\delta}(B,C) \, & = \, f(C)\ts n^{2} \. +
\ts \left[ 2\ts f(BC)\. - \. BC\ts \log\left(1+\frac{1}{C}\right) \right] n^{1+\delta}
\. + \. D\ts n^{2\delta} \, + \, O\bigl(n^{3\delta-1}+n\log n\bigr),
\end{align}
where
$$
D \. := \. f(E) \ts + \ts E\ts \log\frac{(1+C)(BC)^{2}}{C(BC+1)^{2}}
\, - \, \frac{B^{2}C}{2(C+1)} \quad \text{and} \quad E \. :=\, \frac{B^{2}C(C+1)}{(B_{c}-B)(B_{c}+B-2)}\..
$$
			
\smallskip
\nin
$(ii)$ \. For \ts $B>B_{c}$, we have:
\begin{align}
	\log \rT_{n,\delta}(B,C) & \,= \, f(C)\ts n^{2} \. + \. \left[ 2\ts f(B_{c}\ts C) \. - \.
B_{c}\ts C\ts \log\left(1+\frac{1}{C}\right) \right]  n^{1+\delta} \. + \.
O\bigl(n^{2\delta}+n\log n\bigr).
\end{align}
\end{theorem}
	
\smallskip

We prove the theorem in the next section.
Note that for \ts $\frac12 < \de <1$, we obtain three terms in the
asymptotics in the first case, and for \ts $0 < \de <1$, two terms in the second case. It is important to note that the coefficient of the second order term $n^{1+\delta}$ as a function of $B$ attains global maximum at $B=B_{c}$ and decreases in the interval $(B_{c},\infty)$.

In fact, the first two terms of \ts $\log \rT_{n,\delta}(B,C)$ \ts agree with those of \ts $\log \rI_{n,\delta}(B,C)$, see Lemma~\ref{lemma:IH_vol_entropy}. Hence, the independence heuristic predicts that the number \ts $\rT_{n,\delta}(B,C)$ \ts of contingency tables to \textit{decrease} in $B$ when $B>B_{c}$. However, Theorem \ref{thm:main_volume} proves that \ts $\rT_{n,\delta}(B,C)$ \ts remains constant for $B>B_{c}$ (up to the second order) due to the phase transition at $B=B_{c}$, which was `invisible' to the independence heuristic. See $\S$\ref{ss:finrem-over} for further discussion.

\medskip

\subsection{Asymptotic independence}\label{ss:main-reason}
To understand the main theorem, consider the \emph{correlation ratio} in
contingency tables, defined as
$$
\rho(\c,\r) \, := \, \frac{\rT(\r,\c)}{\rI(\r,\c)}\..
$$
%
%
Let us show how \ts $\rho(\c,\r)$ \ts can be interpreted as the asymptotic
independence of rows and columns.  First, recall the following is the argument
essentially in \cite{Good76} (see also \cite{B1}).

Let \ts $\cS_{N}$ \ts be the set of all $m\times n$ tables with total sum
$N=a_1+\ldots + a_m = b_1+\ldots + b_n$.
Let $X$ be a uniformly chosen contingency table from the set $\cS_N$,
and consider the following events that $X$ satisfies the row and column margins:
\begin{align}\label{eq:def_margin_events}
\mathcal{R}_n(\a) \. = \. \bigl\{ \text{$X$ has row margins $\a$} \bigr\}\quad \text{and} \quad
\mathcal{C}_m(\b) \. = \. \bigl\{ \text{$X$ has column margins $\b$} \bigr\}.
\end{align}
By definition
\begin{align}
\P\bigl(\mathcal{R}_n(\r)\cap \mathcal{C}_m(\c)\bigr) \, = \, \frac{\rT(\r,\c)}{|\cS_N|},\qquad
\P\bigl(\mathcal{R}_{n}(\r)\bigr) \, = \, \frac{|\mathcal{R}_{n}(\r)|}{|\cS_N|}, \qquad
\P\bigl(\mathcal{C}_{n}(\c)\bigr) \, = \, \frac{|\mathcal{C}_{n}(\c)|}{|\cS_N|}\..
\end{align}
Since
$$
\bigl|\cS_{N}\bigr| \. = \. \binom{N+mn-1}{mn-1}, \quad \
\bigl|\mathcal{R}_n(\a)\bigr| \. = \.\prod_{i=1}^{m} \. \binom{a_{i}+n-1}{n-1}, \ \quad \,
\bigl|\mathcal{C}_m(\b)\bigr| \. = \. \prod_{j=1}^{n} \. \binom{b_{j}+m-1}{m-1},
$$
using the definition of independence heuristic~\eqref{eq:IH-def} we conclude that
\begin{align}
\frac{\P\bigl(\mathcal{R}_{n}(\a)\cap \mathcal{C}_{m}(\b)\bigr)}{\P\bigl(\mathcal{R}_{n}(\a)\bigr) \, \P\bigl(\mathcal{C}_{m}(\b)\bigr)}
\, = \, \frac{\rT(\r,\c)}{\rI(\r,\c)}\, = \, \rho(\c,\r).
\end{align}
In other words, \ts $\log \rho(\bba,\bbb)$ \ts measures the independence of row
and column margins in random contingency tables.

\medskip

\subsection{Critical correlation exponent}	\label{ss:main-exp}
We now compute the asymptotics of the correlation ratio for the Barvinok
margins~\eqref{eq:def_barvinok_margin}.  As stated in the introduction,
we prove that the ratio also exhibits a second order phase
transition in parameter~$B$.
	
\smallskip
	
\begin{theorem}\label{thm:correlation}
Fix $0< \delta<1$, $B,C>0$, and $n\ge 1$. Let \ts $B_{c}=1+\sqrt{1+1/C}$ \ts
and denote \ts $f(x) := (x+1)\log(x+1) - x\log x$. Then:
\begin{align}\label{eq:thm_correlation_1}
		\qquad
\lim_{n\rightarrow\infty} \. \frac{1}{n^{1+\delta}} \. \log \frac{\rT_{n,\delta}(B,C)}{\rI_{n,\delta}(B,C)} \,  =
		\begin{cases}
		0 &  \text{if \, $B \le B_{c}$} \\
		C(B-B_{c})\log \left(1+\frac{1}{C}\right) \ts - \ts 2\bigl(f(BC) - f(B_{c}C)\bigr) \. > \. 0 &  \text{if \, $B>B_{c}$}
		\end{cases}
\end{align}
\end{theorem}

\smallskip

Corollary~\ref{cor:main} follows immediately from the theorem for $C=1$.
Note that the \emph{critical factor} $n^{1+\de}$ can also be found in the asymptotics
of the total sum~$N$ in the \ts $(n+\lfloor n^{\delta}\rfloor)\times (n+\lfloor n^{\delta}\rfloor)$ \ts contingency tables
with Barvinok margins:
$$N\, = \, \lfloor Cn \rfloor\cdot n \. + \. 2\ts \lfloor BCn \rfloor \cdot \lfloor n^{\delta}\rfloor
\, = \, Cn^{2} \. + \. 2BCn^{1+\delta} \. + \. O(n).
$$

\begin{figure}[hbt]
		\begin{center}
			\includegraphics[height=4.2cm]{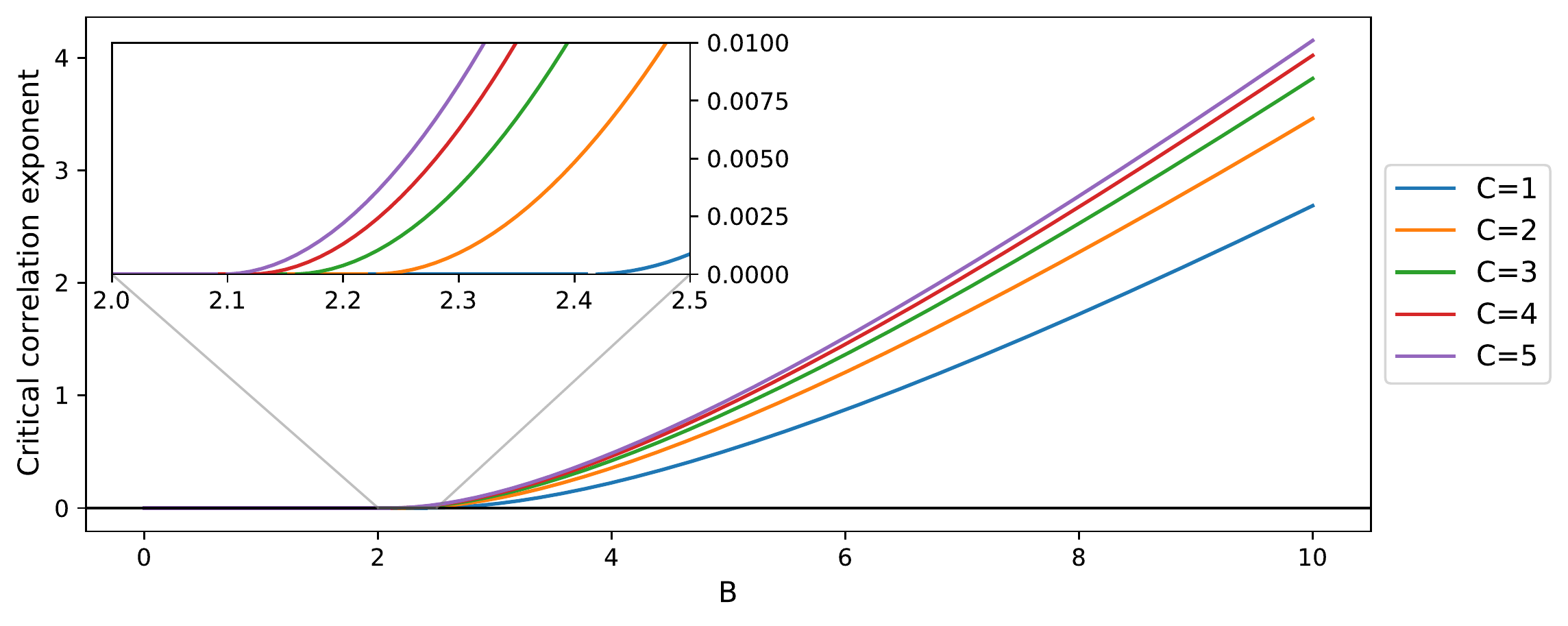}
		\end{center}
		\caption{
			Plot of the critical correlation coefficient in uniform contingency tables with Barvinok's margins with parameters $n,\delta,B$ and $C$. For each $0<\delta<1$ and $C>0$, there exists a second-order phase transition of the critical correlation coefficient in $B$ at critical value $B_{c}=1+\sqrt{1+1/C}$. Below $B_{c}$ the rows and columns are asymptotically independent, but above $B_{c}$, they are asymptotically positively correlated.
		}
		\label{fig:CT}
	\end{figure}

We call the left hand side of \eqref{eq:thm_correlation_1} the \textit{critical correlation exponent}
for the contingency tables $\M(\a,\b)$.  Theorem \ref{thm:correlation} implies that the row and column
margin events $\mathcal{R}_{n}$ and $\mathcal{C}_{n}$ are asymptotically independent for $B<B_{c}$
and asymptotically positively correlated for $B>B_{c}$.
Moreover, it is easy to check that the right hand side of \eqref{eq:thm_correlation_1} as well as its first derivative in $B$ is continuous for $B>0$, but its second derivative is discontinuous at $B=B_{c}$ (see Figure~\ref{fig:CT}). Hence we are uncovering a second-order phase transition in the correlation structure in contingency tables.
	
\medskip

\section{Proof of Theorem~\ref{thm:main_volume}}
	
The proof below relies on the notion of \emph{typical table} introduced in~\cite{B3}
(Definition~\ref{def:typical_table}), and a result in~\cite{DLP19} that was used
prove a probabilistic phase transition for the uniformly sampled contingency table
for the Barvinok margins (Lemma~\ref{lemma:typical_mx}).
	
Let $\mathcal{P}(\r,\c)\subseteq \R_+^{mn}$ be the \textit{transportation polytope}
of \emph{real} nonnegative contingency tables with margins $\r$ and~$\c$, i.e.\
defined by~\eqref{eq:def_contingency_set} over~$\rr_+$.
Clearly, $\M(\r,\c) = \mathcal{P}(\r,\c) \cap \zz^{mn}$.
Next, we define the \textit{typical table} introduced by  Barvinok~\cite{B3}.
	
\smallskip

\begin{definition}[Typical table]\label{def:typical_table}
		Fix margins $\r\in \nn^{m}$ and $\c\in \nn^{n}$.
		Let $\mathcal{P}(\r,\c)\subseteq \R_+^{mn}$ denote the transportation polytope.
		For each $X=(X_{ij})\in \mathcal{P}(\r,\c)$, define
		\begin{align}\label{eq:def_g}
		g(X) \. = \, \sum_{1\le i,j\le n} \. f(X_{i,j})\ts,
		\end{align}
		where the function $f:[0,\infty)\rightarrow [0,\infty)$ is defined by
		\begin{align}\label{eq:def_f}
		f(x) \. = \. (x+1)\log(x+1) - x\log x\ts.
		\end{align}
		The \textit{typical table} $Z\in \mathcal{P}(\r,\c)$ for $\M(\r,\c)$ is defined by
		\begin{align}
		Z \. = \. \argmax_{X\in \mathcal{P}(\r,\c)} \. g(X).
		\end{align}
\end{definition}

Since the function $g$ defined at~\eqref{eq:def_g} is strictly concave, it attains a unique maximizer on the transportation polytope $\mathcal{P}(\r,\c)$ and thus the typical table is well-defined.
	
In \cite[Thm~1.1]{B1}, Barvinok gave the following upper and lower bound on
the number of contingency tables in terms of the typical table
(cf.~\ref{ss:finrem-other}).

\smallskip

\begin{theorem}[\cite{B1}]\label{thm:Barvinok_volume}
Fix margins $\r\in \nn^{m}$ and $\c\in \nn^{n}$. Let $Z=(z_{ij})$ be the typical table for $\M(\r,\c)$.  
Then there exists some absolute constant $\gamma>0$, such that
\begin{align}
g(Z) \. - \. \gamma(m+n) \log N \, \le \, \log \ts \rM(\r,\c) \, \le \, g(Z)\ts,
\end{align}
where \. $N = a_1+\ldots +a_m = b_1+\ldots+b_n$ \. is the total sum of the entries.
\end{theorem}

\smallskip

The following an asymptotic expression of the solution to the optimization problem
for the typical table for $\M_{n,\delta}(B,C)$. A slight modification of the argument shows the following:
	
\smallskip

\begin{lemma}[{\cite[Lem~5.1]{DLP19}}]\label{lemma:typical_mx}
Let $Z=(z_{ij})$ be the typical table for \ts $\M_{n,\delta}(B,C)$, where \ts $0\le \delta<1$.
Let \ts $B_{c}=1+\sqrt{1+1/C}$. Then there exists a constant $\alpha=\alpha(C)$ independent of~$B$,
such that:
\begin{description}
\item[(i)] \, If \. $B \le B_{c}$\ts, then:
\begin{align}
\bigl|z_{n+1,n+1} - C\bigr| \, & \le \, BC \ts n^{\delta-1} \\
\bigl|z_{1,n+1} - BC\bigr| \, &\le \, \frac{\alpha}{B_{c}-B} \ts n^{\delta-1} \\
\left| z_{11} - \frac{B^{2}C(C+1)}{ (B_{c}-B)(B_{c}+B-2) }\right| \, & \le \, \frac{\alpha}{B_{c}-B} \ts n^{\delta-1}.
\end{align}
			
\item[(ii)]\,  If \. $B>B_{c}$\ts, then:
\begin{align}
\bigl|z_{n+1,n+1} - C\bigr| \, & \le \, B_{c}C \ts n^{\delta-1} \\
\bigl|z_{1,n+1} - BC\bigr| \, &\le \, \frac{\alpha}{B-B_{c}} \ts n^{\delta-1} \\
\left| n^{\delta-1}z_{11} - C(B-B_{c}) \right| \, & \le \, \frac{\alpha}{B-B_{c}} \. n^{\delta-1}.
		\end{align}
		\end{description}
	\end{lemma}

\smallskip

We use the lemma to prove the following result.
	
\smallskip
	
\begin{prop}\label{prop:entropy_estimate}
Let $Z=(z_{ij})$ be the typical table for $\T(\r,\c)$, where $0\le \delta<1$. Let $B_{c}=1+\sqrt{1+1/C}$. Let $f,g$ be the functions defined at \eqref{eq:def_f} and \eqref{eq:def_g}. Then the following hold:
\begin{description}
			\item[(i)] \. If \ts $B<B_{c}$, then for all \ts $n\ge 1$,
			\begin{align}
			g(Z) \, & = \, f(C)\ts n^{2} \. + \. \left[ 2f(BC) \. - \. BC\log\left(1+\frac{1}{C}\right) \right] n^{1+\delta} \\
			&\qquad + \. \left[ f(z_{11}) \. + \. z_{11}^{*}\log\left( \frac{(1+C)(BC)^{2}}{C(BC+1)^{2}}\right)
            \. - \. \frac{B^{2}C}{2(C+1)}  \right] n^{2\delta} + O(n^{3\delta-1}) + O(n),
			\end{align}
			where \ts $z_{11}^{*}=B^{2}C(C+1)/(B_{c}-B)(B_{c}+B-2)$.
			
			\item[(ii)] \. If \ts $B>B_{c}$, then for all \ts $n\ge 1$,
			\begin{align}
			g(Z) \, & = \, f(C)\ts n^{2} \. + \. \left[ 2f(B_{c}C) \. - \.
        B_{c}C\log\left(1+\frac{1}{C}\right) \right] n^{1+\delta} \. + \. O\bigl(n^{2\delta}\bigr) + O(n).
			\end{align}
\end{description}
\end{prop}

\smallskip

	\begin{proof}
		First recall that due to the symmetry, the entries \ts $z_{11}$, $z_{1,n+1}$ \ts
and \ts $z_{n+1,n+1}$ \ts of the typical table~$Z$ satisfy the following margin condition:
\begin{align}\label{eq:z_ij_reduced}
\begin{cases}
		(\lfloor n^{\delta} \rfloor/n) \ts z_{11} + z_{1,n+1} \. = \. (\lfloor BCn\rfloor /n) = BC + O(n^{-1})\ts, \\
		(\lfloor n^{\delta} \rfloor/n) \ts z_{1,n+1} + z_{n+1,n+1} \. = \. (\lfloor Cn\rfloor /n) = C + O(n^{-1})\ts.
\end{cases}
\end{align}
In accordance to Lemma~\ref{lemma:typical_mx}, define a block table $Z^{*}=(z_{ij}^{*})$ by $z_{n+1,n+1}^{*}=C$, $z_{1,n+1}=BC$, and $z_{11}^{*}=B^{2}C(C+1)/(B_{c}-B)(B_{c}+B-2)$. Combining with Lemma~\ref{lemma:typical_mx}, for $B<B_{c}$, we obtain
		\begin{align}
		C-	z_{n+1,n+1} &= z_{1,n+1}(\lfloor n^{\delta}\rfloor/n) + O(n^{-1}) \\
		&= BC(\lfloor n^{\delta}\rfloor/n) - (BC-z_{1,n+1})(\lfloor n^{\delta}\rfloor/n) + O(n^{-1})\\
		&= BC(\lfloor n^{\delta}\rfloor/n)  -  z_{11}^{*}(\lfloor n^{\delta}\rfloor/n)^{2} - (z_{11}-z_{11}^{*})(\lfloor n^{\delta}\rfloor/n)^{2} + O(n^{-1}) \\
		&=  BC(\lfloor n^{\delta}\rfloor/n)  -  z_{11}^{*}(\lfloor n^{\delta}\rfloor/n)^{2} +  O(n^{3\delta-3}) + O(n^{-1}),
		\end{align}
		and also
		\begin{align}
		BC-z_{1,n+1} &=  z_{11}^{*}(\lfloor n^{\delta}\rfloor/n) + (z_{11}-z_{11}^{*})(\lfloor n^{\delta}\rfloor/n) + O(n^{-1})\\
		&= z_{11}^{*}(\lfloor n^{\delta}\rfloor/n) +  O(n^{2\delta-2}) + O(n^{-1}).
		\end{align}
		Similarly, for $B>B_{c}$,
		\begin{align}
		C-z_{n+1,n+1} &= B_{c}C(\lfloor n^{\delta}\rfloor/n) + (z_{1,n+1}-B_{c}C)(\lfloor n^{\delta}\rfloor/n) + O(n^{-1})\\
		& = B_{c}C(\lfloor n^{\delta}\rfloor/n)  + \left[ C(B-B_{c}) - (\lfloor n^{\delta}\rfloor/n)z_{11} \right](\lfloor n^{\delta}\rfloor/n) + O(n^{-1})\\
		&=  B_{c}C(\lfloor n^{\delta}\rfloor/n) + O(n^{2\delta-2}) + O(n^{-1}).
		\end{align}

		Now suppose $B \le B_{c}$. We use the following Taylor expansion of $f$:
		\begin{align}
		f(x) \. = \. f(y) + (x-y) \log(1+x^{-1}) + \frac{(x-y)^{2}}{2y(y+1)} + O(|x-y|^{3}),
		\end{align}
		where the constant in $O(\cdot)$ above is bounded when $x,y>0$ are remain bounded. Then observe that
		\begin{align}
		f(z_{n+1,n+1}) &= f(C) + (z_{n+1,n+1}-C)\log\left(1+\frac{1}{C}\right) + \frac{(z_{n+1,n+1}-C)^{2}}{2C(C+1)} + O(n^{3\delta-3}) + O(n^{-1}), \\
		&=f(C)-BC (\lfloor n^{\delta} \rfloor/n) + z_{11}^{*}(\lfloor n^{\delta} \rfloor/n)^{2} + O(n^{3\delta-3}) + O(n^{-1}), \\
		f(z_{1,n+1}) &=f(BC) + (z_{1,n+1}-BC) \log(1+1/BC) + O(n^{2\delta-2}) + O(n^{-1}) \\
		&= f(BC) - z_{11}^{*} (\lfloor n^{\delta} \rfloor/n) + O(n^{2\delta-2}) + O(n^{-1}), \\
		f(z_{11}) &= f(z_{11}^{*}) + O(n^{\delta-1}) + O(n^{-1}).
\end{align}
Noting that \ts $g(Z)=n^{2}f(z_{n+1,n+1}) + 2n\lfloor n^{\delta} \rfloor f(z_{1,n+1}) + \lfloor n^{\delta} \rfloor^{2}f(z_{11})$, a straightforward computation shows \textbf{(i)}. Next, suppose $B>B_{c}$. By a similar argument, we have
\begin{align}
		f(z_{n+1,n+1}) &= f(C) - B_{c}Cn^{\delta-1} \log\left(1+\frac{1}{C}\right) + O(n^{2\delta-2}) + O(n^{-1}) ,  \\
		f(z_{1,n+1}) &=f(B_{c}C) + O(n^{\delta-1}),\qquad f(z_{11}) = f(z_{11}^{*}) + O(n^{\delta-1}) + O^{(n^{-1})}.
\end{align}
Then \textbf{(ii)} follows from here.
\end{proof}

\smallskip
		
	
\begin{proof}[\textbf{Proof of Theorem~\ref{thm:main_volume}}]
Suppose \ts $0< \delta < 1$. Note that \ts $N=Cn^{2}+BCn^{1+\delta}+O(n)$,
where $N$ denotes the total sum of entries in a contingency table in \ts
$\M_{n,\delta}(B,C)$. Combining with Theorem \ref{thm:Barvinok_volume}, we have:
\begin{align}
\bigl| \log \rT(\a,\b) - g(Z) \bigr| \, \le \, \gamma' n\log n\ts,
\end{align}
for all \ts $n\ge 1$, for some absolute constant \ts $\gamma'>0$.
Now the theorem follows from Proposition~\eqref{prop:entropy_estimate}.
\end{proof}
		
\bigskip
		
\section{Proof of Theorem~\ref{thm:correlation}}

We start with the following lemma.
	
\smallskip

\begin{lemma}\label{lemma:IH_vol_entropy}
We have:
\begin{align}
\log \. \rI_{n,\delta}(B,C) \, & = \,
f(C)n^{2} + \left[ f(BC) - BC\log\left(1+\frac{1}{C}\right) \right]n^{1+\delta}  \\
&\qquad + \left[ 2\log(BC+1) - \log (C+1) + \frac{(2-4B+B^{2})C}{2(1+C)} \right]n^{2\delta}  +O(n^{3\delta-1}+n\log n).
\end{align}
\end{lemma}
	
\smallskip

The proof of Lemma~\ref{lemma:IH_vol_entropy} involves a straightforward
computation of expanding the right hand side of \eqref{eq:IH-def} under
Barvinok's margins \eqref{eq:def_barvinok_margin}.
Details are given in the next Section~\eqref{section:pf_IH_estimate}.
	
\smallskip

\begin{proof}[\textbf{Proof of Theorem \ref{thm:correlation}}]
		By Lemma~\ref{lemma:IH_vol_entropy},
		\begin{align}\label{eq:IH_vol_formula_pf_in_thm}
		\log \, \rI_{n,\delta}(B,C)&= f(C)n^{2} +  \left[ 2f(BC) - BC\log\left(1+\frac{1}{C}\right)\right]n^{1+\delta} + O(n^{3\delta}+n\log n).
		\end{align}
		Suppose $B<B_{c}$. Recall that by Theorem \ref{thm:main_volume}~\textrm{(i)},
		\begin{align}
		\log \rT_{n,\delta}(B,C) \, &= \, f(C)n^{2} \. + \.
        \left[ 2f(BC) - BC\log\left(1+\frac{1}{C}\right) \right] n^{1+\delta} \. + \. O(n^{2\delta}+n\log n),
		\end{align}
		Hence for \ts $0\le \delta <1$, we have \ts $1+\delta>2\delta$, so we obtain:
		\begin{align}
		\lim_{n\rightarrow\infty} \. \frac{1}{n^{1+\delta}}\. \log \frac{ \rT(\r, \c)}{ \rI(\r,\c)} \, = \, 0.
		\end{align}
		On the other hand, suppose \ts $B>B_{c}$. Then by Theorem~\ref{thm:main_volume}~\textrm{(ii)},
		\begin{align}
		\log \rT_{n,\delta}(B,C) \, =  \, f(C)n^{2} + \left[ 2f(B_{c}C)  - B_{c}C\log\left(1+\frac{1}{C}\right) \right] n^{1+\delta} + O(n^{2\delta}+n\log n).
		\end{align}
		Hence by \eqref{eq:IH_vol_formula_pf_in_thm}, for \ts $0<\delta <1$, we have:
		\begin{align}
		\lim_{n\rightarrow\infty} \. \frac{1}{n^{1+\delta}} \.
        \log \. \frac{ \rT_{n,\delta}(B,C)}{\rI_{n,\delta}(B,C)} \,
        &= \, 2f(B_{c}C) - 2f(BC)  + C(B-B_{c})\log\left(1+\frac{1}{C}\right).
		\end{align}
		Let us denote the right hand side of the above equation as $\lambda(B)$. Then
		\begin{align}
		\frac{\partial \lambda}{\partial B}  \, &= \, C\log\left(1+\frac{1}{C}\right) \. - \. 2C\log\left(1+\frac{1}{BC}\right),
		\end{align}
		and it is easy to see $\frac{\partial \lambda}{\partial B}>0$ if and only if $B>B_{c}$ and the derivative at $B=B_{c}$ equals zero. Hence for each fixed $C$, the function $\lambda$ is a strictly increasing on $[B_{c},\infty)$, and has minimum at $B=B_{c}$.  Note also that $\la(B)=0$.  This shows $\lambda(B)>0$ for all $B>B_{c}$.
	\end{proof}
	
We remark that the second derivative in $B$ of the limiting expression in \eqref{eq:thm_correlation_1} for $B>B_{c}$ is $2C/(B(BC+1))>0$ which is also strictly positive at $B=B_{c}$. Hence the phase transition given in Theorem~\ref{thm:correlation} is indeed of second order.

\bigskip

\section{Proof of Lemma~\ref{lemma:IH_vol_entropy}}
\label{section:pf_IH_estimate}

%

	We first compute $\log \rI(\a,\b)$ for the general $m\times n$ tables with total sum~$N$. By Stirling's approximation,
	\begin{align}
	\log \binom{a + b}{a} \, = \, (a+b) \log (a+b) - a\log a -b\log b + O(\log (a+b)).
	\end{align}
%
Then we have:
\begin{equation}\label{eq:IH_volume_expansion}
\begin{aligned}
	& \log \rI(\a,\b) \, = \,  \sum_{i=1}^{\wb{m}} \. (r_{i}+\wb{n})\log (r_{i}+\wb{n}) \, + \,\sum_{j=1}^{\wb{n}} \. (c_{j}+\wb{m})\log (c_{j}+\wb{m}) \, - \, \sum_{i=1}^{\wb{m}} \. r_{i}\log r_{i} \\
	&  - \, \sum_{j=1}^{\wb{n}} \. c_{j}\log c_{j} \, - \, (N+\wb{m}\cdot \wb{n})\log (N+\wb{m}\cdot \wb{n})  + N\log N + O\bigl((\wb{n}+\wb{m})\log(N+\wb{m}\cdot \wb{n})\bigr).
\end{aligned}
\end{equation}
	
	Now assume Barvinok margins \eqref{eq:def_barvinok_margin}. Denote \ts $\phi(x):=x\log x$. We use \eqref{eq:IH_volume_expansion}
with \ts $\wb{m} \gets n+n^{\delta}$, \ts $\wb{n}\gets n+n^{\delta}$, and \ts $N \gets Cn^{2}+BCn^{1+\delta}$.  We have:
	\begin{align}
	\log \rI_{n,\delta}(B,C) \, &= \, 2n\left[ \phi\left( (C+1)n+n^{\delta} \right) - \phi(Cn)\right] + 2n^{\delta}\left[ \phi\left( (BC+1)n+n^{\delta}\right) -\phi\left( BCn \right) \right]   \\
	&\qquad -\phi\left(N+mn\right) + \phi\left(N\right) + O(n\log n).
	\end{align}
	Using Taylor expansion, we have:
	\begin{align}
	& \log \left( (C+1)n + n^{\delta} \right) 
	 \, = \,  \log n + \log (C+1) \. + \. \frac{n^{\delta-1}}{C+1} \. - \. \frac{n^{2\delta-2}}{2(C+1)^{2}} \. + \. O(n^{3\delta-3}), \\
	& \log \left( Cn^{2} + BCn^{1+\delta}\right) 
 \, = \, \log Cn^{2} + Bn^{\delta-1}- \frac{B^{2}n^{2\delta-2}}{2} + O(n^{3\delta-3}),\\
	& \log \left(N+mn\right) 
	\, = \, \log (C+1)n^{2} + \frac{BC+2}{C+1}n^{\delta-1} + \frac{(-1+C-2BC-B^{2}C^{2})n^{2\delta-2}}{(1+C)^{2}} + O(n^{3\delta-3}).
	\end{align}
	Hence we get
	\begin{align}
	& \log \rI_{n,\delta}(B,C) \,  = \, 2n\left((C+1)n+n^{\delta}\right)\left[\log (C+1)n + \frac{n^{\delta-1}}{C+1} - \frac{n^{2\delta-2}}{2(C+1)^{2}}\right]  - 2Cn^{2}\log Cn\\
	& \quad\quad + 2n^{\delta}\left( (BC+1)n+n^{\delta}\right)\left[\log (BC+1)n + \frac{n^{\delta-1}}{BC+1} - \frac{n^{2\delta-2}}{2(BC+1)^{2}} \right] -2n^{1+\delta} BC\log BCn \\
	&\quad\quad- \left((C+1)n^{2} + (BC+2)n^{1+\delta} + n^{2\delta}\right) \. \times \\
	&\quad\quad\quad\qquad \times \. \left[ \log (C+1)n^{2} +  \frac{BC+2}{C+1}n^{\delta-1} +\left[ -1+C-2BC-B^{2}C^{2} \right] \frac{n^{2\delta-2}}{(1+C)^{2}}  \right]\\
	&\quad\quad  + \left( Cn^{2} + BCn^{1+\delta}\right)\left[ \log Cn^{2} + Bn^{\delta-1} - \frac{B^{2}n^{2\delta-2}}{2} \right]+O(n^{3\delta-1}+n\log n)\\
	&= \, f(C)n^{2} + \left[ f(BC) - BC\log\left(1+\frac{1}{C}\right) \right]n^{1+\delta}  \\
	&\qquad + \, \left[ 2\log(BC+1) - \log (C+1) + \frac{(2-4B+B^{2})C}{2(1+C)} \right]n^{2\delta} +O(n^{3\delta-1}+n\log n).
	\end{align}
	This completes the proof. \ $\sq$

\bigskip

\section{Final remarks}\label{sec:finrem}

\subsection{}\label{ss:finrem-hist}
The story behind the \emph{independence heuristic}~\eqref{eq:IH}
and~\eqref{eq:IH-def} is rather interesting.
This approximation was given implicitly by Good in \cite[p.~100]{Good-book},
and later stated formally in~\cite{Good,Good76}.
Good writes that ``the conjecture appears to be confirm''
by his calculations \cite[p.~1166]{Good76}, but later admits
he is ``leaving aside finer points of rigor'' (ibid, p.~1184).

Because of the small constant similar to the $\sqrt{e}$ in~\eqref{eq:CM},
there was an effort to ``improve'' upon~\eqref{eq:IH}.  Unfortunately,
from the asymptotic point of view, many such heuristics behave poorly.
Notably, Diaconis and Efron \cite[$(3.14)$]{DE}, see also~\cite[$(7.2)$]{DG},
propose another heuristic estimate.  For linear margins $a_i,b_j=\Theta(n)$ and
$m=n$, our calculation shows that this formula gives \ts $\log \rTab \approx \Theta(n^2\log n)$,
thus implying the wrong leading term of the asymptotics.

A number of papers tested the independence heuristic numerically,
see e.g.~\cite{DG,GC}.  The results are nothing short of remarkable, showing
that~\eqref{eq:IH} holds up to a small constant.
For example, for the $4\times 4$ case with $N=592$ introduced
in~\cite{DE}, we have \ts $\rM(\r,\c)=1.226 \times 10^{15}$~\cite{DG},
while \ts $\rI(\r,\c)=1.211 \times 10^{15}$~\cite{B1}.  The asymptotic
analysis shows that this level of agreement is purely coincidental.
In fact, for a much larger $4\times 4$ case with $N=65159458$ computed
in~\cite{DL1,DL2}, we have \ts $\rM(\r,\c)=4.3 \times 10^{61}$ \ts
vs.\ \ts $\rI(\r,\c)=3.7 \times 10^{61}$.  This and other numerical
estimates are collected in~\cite[$\S$11.4]{BLP}.

\subsection{} \label{ss:finrem-other}
There are several better lower and upper bounds known for the number \ts $\rTab$ \ts
of contingency tables, see~\cite{BLP} for a recent overview.  Notably, the upper bound
in~\cite{Sha} and the most recent lower bound in~\cite{BLP} give improvement over
Theorem~\ref{thm:Barvinok_volume} in the lower order terms.  For the linear margins
these improvement are of the order $n^{cn}$. Thus, they give an improvement
in the second order term in the Main Theorem~\ref{ss:main-thm}~$(i)$ when $\de<\frac12$.
It would be interesting to further explore these bounds.

\subsection{} \label{ss:finrem-general}
Consider $m\times n$ contingency tables with general margins $\r$ and $\c$ and total sum~$N$.
Let $W=(w_{ij})$ be as $w_{ij}=r_{i}c_{j}/N$ for $1\le i \le m$, $1\le j \le n$.
Let $g$ be defined in \eqref{eq:def_g}. In \cite[$\S$2]{B1}, Barvinok showed that
\ts $g(W) - \log \rI(\a,\b) \ge 0$. This implies that
\begin{align}\label{eq:barvinok_entropy_ineq}
		\log \, \frac{\rT(\r,\c)}{\rI(\r,\c)} \, \ge \, - \. \gamma \. (m+n) \ts \log N.
\end{align}
This shows that the row and column margin events $\mathcal{R}$ and $\mathcal{C}$
defined by~\eqref{eq:def_margin_events}, have asymptotically \textit{nonnegative correlation} at a scale where the right hand side of \eqref{eq:barvinok_entropy_ineq} vanishes. For instance, for $m=n$ and linear margins $a_i, b_j = \Theta(n)$, we have $\gamma (m+n)\log N = O(n\log n)$. Then:
\begin{align}\label{eq:barvinok_entropy_ineq2}
    \limsup_{n\rightarrow \infty}\. \frac{1}{n^{\delta+1}}\. \log \, \frac{\rT(\r,\c)}{\rI(\r,\c)} \, \ge \, 0,
\end{align}
for every $\delta>0$.
	

\subsection{} \label{ss:finrem-over}
To understand the discussion of the Main Theorem~\ref{thm:main_volume} at the end of~$\S$\ref{ss:main-thm},
consider an extreme case of the margins \ts $\a=\b=(Bn^2,n,\ldots,n)\in \rr^{n+1}$. For every fixed~$n$,
when $B>1$ is large enough, the number $\rTab$ stabilizes as the corner entry \ts $x_{11}$ \ts
is forced to absorb bulk of the total sum \ts $N = (B+1)n^2$.  Meanwhile, the independence
heuristic approximation \ts $\rI(\a,\b)$ \ts maximizes at a certain constant~$B$ \ts
and then decreases exponentially, eventually becoming~$<1$.  In other words, our
Main Theorem~\ref{thm:main_volume} implies a lower order version of the same phenomenon.

Note also a different but related phenomenon of the lower bound in Barvinok's
Theorem~\ref{thm:Barvinok_volume}, which works well asymptotically for large~$n$, but
for small $n$ and large marginals gives lower bonds which are~$<1$, see~\cite[$\S10$]{BLP}.

\subsection{} \label{ss:finrem-bin}  
In~\cite{B2}, Barvinok gave the analogue of Theorem~\ref{thm:Barvinok_volume} for 
\emph{binary} (0-1) \emph{contingency tables}.   Most recently, Wu~\cite{Wu} 
investigated the limiting distribution of the entries in binary contingency tables, 
and proved the analogue of our Lemma~\ref{lemma:typical_mx}.  It would be interesting
to find the analogue of our Main Theorem~\ref{thm:main_volume} in this setting.

\vskip.3cm

{\small
\subsection*{Acknowledgements}
We are grateful to Sasha Barvinok, Sam Dittmer, Jonathan Leake and Greta Panova for numerous
helpful discussions, and to Robin Pemantle for telling us about paper~\cite{Wu}. 
Both authors were partially supported by the NSF.
}

\vskip1.5cm



{\footnotesize
		

\begin{thebibliography}{BLSY10}

\bibitem[B1]{B1}
A.~Barvinok,
Asymptotic estimates for the number of contingency tables, integer flows, and volumes of transportation polytopes,
\emph{Internat.\ Math.\ Res.\ Notices}~\textbf{2009} (2009), 348--385.

\bibitem[B2]{B2}
A.~Barvinok,
On the number of matrices and a random matrix with prescribed
row and column sums and 0-1 entries,
\emph{Adv.\ Math.}~\textbf{224} (2010), 316--339.

\bibitem[B3]{B3}
A.~Barvinok, What does a random contingency table look like?,
\emph{Comb.\ Probab.\ Comp.}~\textbf{19} (2010), 517--539.

\bibitem[BH]{BH}
A.~Barvinok and J.~A.~Hartigan,
An asymptotic formula for the number of non-negative integer matrices with prescribed row and column sums,
\emph{Trans.\ AMS}~\textbf{364} (2012), 4323--4368.
			
\bibitem[BBK]{BBK}
A.~B\'{e}k\'{e}ssy, P.~B\'{e}k\'{e}ssy and J.~Koml\'{o}s,
Asymptotic enumeration of regular matrices,
\emph{Studia Sci.\ Math.\ Hungar.}~\textbf{7} (1972), 343--353.
			
\bibitem[BLP]{BLP}
P.~Br\"and\'en, J.~Leake and I.~Pak,
Lower bounds for contingency tables via Lorentzian polynomials,
preprint (2020), 28 pp.; \ts {\tt arXiv:2008.05907}.

\bibitem[CM]{CM}
E.~R.~Canfield and B.~D.~McKay,
Asymptotic enumeration of integer matrices with large equal row and column sums,
\emph{Combinatorica}~\textbf{30} (2010), 655--680.
			
			
\bibitem[D1]{DL1}
J.~A.~De~Loera, Counting and Estimating Lattice Points:
Tools from Algebra, Analysis, Convexity,
and Probability, \emph{Optima}~\textbf{81} (2009), 1--9; available at
\ts \url{http://www.mathopt.org/Optima-Issues/optima81.pdf}

\bibitem[D2]{DL2}
J.~A.~De~Loera, Details on experiments (counting and estimating
  lattice points), an appendix to~\cite{DL1}, ibid., 17-22;  available at
\ts \url{http://www.mathopt.org/Optima-Issues/optima81-app.pdf}

			
\bibitem[DE]{DE}
P.~Diaconis and B.~Efron,  Testing for independence in a two-way table:
new interpretations of the chi-square statistic,
\emph{Ann.\ Stat.}~\textbf{13} (1985), 845--913.
			
\bibitem[DG]{DG}
P.~Diaconis and A.~Gangolli,
Rectangular arrays with fixed margins,
\emph{Disc.\ Prob.\ Alg.}~\textbf{72} (1995), 15--41.
			
			
\bibitem[DLP]{DLP19}
S.~Dittmer, H.~Lyu and I.~Pak,
Phase transition in random contingency tables with non-uniform margins,
to appear in \emph{Trans.\ AMS}; \ts {\tt arXiv:1903.08743}.

\bibitem[G1]{Good-book}
I.~J.~Good, \emph{Probability and the Weighing of Evidence},
Hafner, New York, 1950, 119~pp.

\bibitem[G2]{Good}
I.~J.~Good, On the application of symmetric Dirichlet distributions
and their mixtures to contingency tables,
\emph{Annals Math.\ Stat.}~\textbf{34} (1963), 911--934.
			
\bibitem[G3]{Good76}
I.~J.~Good, Maximum entropy for hypothesis formulation, especially for
multidimensional contingency tables,
\emph{Annals.\ Stat.}~\textbf{4} (1976), 1159--1189.
			
\bibitem[GC]{GC}
I.~J.~Good and J.~F.~Crook,
The enumeration of arrays and a generalization related to contingency tables,
\emph{Discrete Math.}~\textbf{19} (1977), 23--45.
		
\bibitem[GM]{GM}
C.~Greenhill and B.~D.~McKay,
Asymptotic enumeration of sparse nonnegative integer matrices with specified
row and column sums, \emph{Adv.\ Appl.\ Math.}~\textbf{41} (2008), 459--481.
			
\bibitem[LW]{LW}
A.~Liebenau and N.~Wormald,
Asymptotic enumeration of digraphs and bipartite graphs by degree sequence,
preprint (2020), 25~pp.; {\tt arXiv:2006.15797}.

\bibitem[Sha]{Sha}
A.~Shapiro, Bounds on the number of integer points in a polytope
via concentration estimates, preprint (2010), 23~pp.; {\tt arXiv:1011.6252}.

\bibitem[Wor]{W2}
N.~Wormald, Asymptotic enumeration of graphs with given degree sequence, in
\emph{Proc.\ ICM Rio de Janeiro}, Vol.~3, 2018, 3229--3248.			

\bibitem[Wu]{Wu}
Da Wu, 
On limiting distribution of a certain class of random binary contingency tables, 
preprint (2020), 11~pp.; {\tt arXiv:2002.12559}. 

\end{thebibliography}		
}
\end{document}